\documentclass[submission,copyright,creativecommons]{eptcs}
\providecommand{\event}{QPL 2015}

\usepackage{subcaption}
\usepackage{mathtools}
\usepackage{youngtab}
\usepackage{multicol}
\usepackage{mathrsfs}
\usepackage{hyperref}
\usepackage{stmaryrd}
\usepackage{amsmath}
\usepackage{amssymb}
\usepackage{amsthm}
\usepackage{braket}

\newcommand\numberthis{\addtocounter{equation}{1}\tag{\theequation}}

\newcommand{\csubstack}[1]{\mathclap{\substack{#1}}}

\newcommand{\Tr}{\mathop{\mathrm{Tr}}}
\newcommand{\after}{\mathbin{\circ}}
\newcommand{\IM}{\mathop{\mathrm{Im}}}
\newcommand{\GL}{\mathop{\mathrm{GL}}}

\newcommand{\id}{\mathrm{id}}
\newcommand{\op}{\mathrm{op}}
\newcommand{\cat}[1]{\mathbf{#1}}
\newcommand{\starpu}{\cat{Star}_\mathrm{cPU}}
\newcommand{\repsn}{\cat{Rep}(S_n)}
\newcommand{\breppsn}{\cat{BRep}({\textstyle\prod}_n S_n)}
\newcommand{\repG}{\cat{Rep}(G)}
\newcommand{\cdn}{(\mathbb{C}^d)^{{\otimes n}}}
\newcommand{\C}{\mathbb{C}}
\newcommand{\N}{\mathbb{N}}
\newcommand{\Hil}{\mathscr{H}}
\newcommand{\Sil}{\mathscr{S}}
\newcommand{\Ail}{\mathscr{A}}

\newcommand{\Sets}{\cat{Set}}

\theoremstyle{definition}

\theoremstyle{remark}

\title{Unordered Tuples in Quantum Computation}

\author{%
Robert Furber
\institute{Radboud Universiteit Nijmegen}
\email{rfurber@cs.ru.nl} \and
Bas Westerbaan
\institute{Radboud Universiteit Nijmegen}
\email{bwesterb@cs.ru.nl}}

\def\authorrunning{Furber \& Westerbaan}
\def\titlerunning{Unordered Tuples in Quantum Computation}

\begin{document}

\maketitle

\begin{abstract}
    It is well known that
    the~C$^*$-algebra of an ordered pair of qubits is $M_2 \otimes M_2$.
        What about \emph{unordered} pairs?
    We show in detail that~$M_3 \oplus \mathbb{C}$
        is the C$^*$-algebra of an \emph{unordered} pair of qubits.
    Then we use Schur-Weyl duality
    to characterize the~C$^*$-algebra of
        an unordered~$n$-tuple of~$d$-level quantum systems.
    Using some further elementary representation theory and number theory,
        we characterize the quantum cycles.
    We finish with a characterization of the von Neumann algebra
        for unordered words.
\end{abstract}

Finite dimensional quantum computation is naturally viewed as
occurring in the category of finite dimensional C$^*$-algebras
together with completely positive unital maps, in the opposite
of their usual direction.
The C$^*$-algebras are the types (systems).  For instance:
\begin{center}\begin{tabular}{lllll}
        a single qubit & $M_2$ &\qquad&
        an ordered pair of qutrits & $M_3 \otimes M_3$ \\
        a single qutrit & $M_3$ &&
        a bit & $\mathbb{C}^2$ \\
        a qubit or
            a qutrit & $M_2 \oplus M_3$ &&
        a trit or a qubit & $\mathbb{C}^3 \oplus M_2$
\end{tabular}\end{center}
More generally, writing~$\llbracket \textrm{type} \rrbracket$
for the C$^*$-algebra for the type, we have:
\begin{align*}
    \llbracket \textrm{$d$-level quantum system} \rrbracket
            & = M_d &
    \llbracket \textrm{$d$-level classical system} \rrbracket
            & = \mathbb{C}^d \\
    \llbracket \textrm{(ordered) pair of $t$ and~$s$} \rrbracket
            & = \llbracket t \rrbracket \otimes \llbracket s \rrbracket &
    \llbracket \textrm{$t$ (classical) or~$s$} \rrbracket
            & = \llbracket t \rrbracket \oplus \llbracket s \rrbracket.
\end{align*}
The completely positive unital maps
are the programs (operations) in the opposite direction.
For example:
\begin{multicols}{2}
\begin{enumerate}
    \item \emph{Measure a qubit in the standard basis} \\
        $m \colon M_2 \leftarrow \C^2$ \qquad (qubit $\to$ bit) \\
        $(\lambda, \mu) \mapsto \lambda \ket{0}\bra{0} + \mu \ket{1}\bra{1}$
    \item \emph{Apply Hadamard gate to a qubit} \\
        $h \colon M_2 \leftarrow M_2 $ \qquad
                (qubit $\to$ qubit) \\
    $a \mapsto H^\dagger a H$,
        where~$H = \frac{1}{\sqrt{2}}\left(\begin{smallmatrix}
                1 & 1 \\ 1 & -1
        \end{smallmatrix}\right)$
    \item \emph{Initialize a qutrit as 0} \\
        $i\colon \C \leftarrow M_3$ \qquad (empty $\to$ qutrit) \\
        $a \mapsto \bra{0}a\ket{0}$
    \item \emph{Forget about a qubit} \\
        $d\colon M_2 \leftarrow \C$ \qquad (qubit $\to$ empty) \\
        $\lambda \mapsto \lambda 1$
\end{enumerate}
\end{multicols}
These basic quantum types are well known, but what about an unordered pair of qubits?
An unordered pair of bits is simply a
trit ($00$, $01=10$ or $11$). 
However, we will see an unordered pair of qubits is not
a qutrit, but rather its C$^*$-algebra is~$M_3 \oplus \C$.

In Section~\ref{s:upair-qubits},
we prove this in detail
to get a feel for this surprising result.
Then in Section~\ref{s:utuples},
we characterize the C$^*$-algebras
of unordered~$n$-tuples of~$d$-level quantum systems
using Schur-Weyl duality, which dates back to the early 20th century.
    Applying some elementary representation theory,
        we characterize the~C$^*$-algebra
        for qubit~$3$-cycles in Section~\ref{s:3cycle}.
        Then, using some number theory,
        we characterize arbitrary quantum cycles
        in Section~\ref{s:cycles}.
    We finish with a characterization of the von Neumann algebra
        for the quantum unordered words in Section~\ref{s:uwords}.

Unordered quantum types have been considered before.  For instance
in~\cite{bew143} they are used to give denotational semantics to a
quantum lambda calculus.\footnote{The
type~$\llbracket \mathbf{qubit} \rrbracket^{\odot 2} $
from \cite[Example 23]{bew143} corresponds to
an unordered pair of qubits and thus has as~C$^*$-algebra
$M_3\oplus \C$.}
A concrete description, however,
has to our knowledge not been published before.

At the end of the paper we will have demonstrated the following.

        \begin{center}\fbox{\begin{minipage}{0.8\textwidth}
        \begin{center}\begin{tabular}{ll}
            System & Algebra \\ \hline
            unordered pair of qubits &
                $M_3 \oplus \C$ \\
            \begin{tabular}[x]{@{}l@{}}
                unordered~$n$-tuple \\
                    \quad of~$d$-level quantum systems
            \end{tabular}&
                $\displaystyle \bigoplus_{\lambda \in Y_n} M_{m_\lambda}$ \\
            words of qubits &
                $B(\ell^2)$ \\
            \begin{tabular}[x]{@{}l@{}}
                unordered words \\
                    \quad of~$d$-level quantum systems
            \end{tabular}&
                $\displaystyle B(\ell^2) \oplus
                    \bigoplus_{\lambda \in Y^*} M_{m_\lambda}$\\
            $3$-cycle of qubits &
                    $M_4 \oplus M_2 \oplus M_2$ \\
            \begin{tabular}[x]{@{}l@{}}
                $n$-cycle \\
                    \quad of~$d$-level quantum systems
            \end{tabular}&
                $\displaystyle \bigoplus_{0 \leq k < n} M_{c_k}$
        \end{tabular}\end{center}
            \begin{align*}
        Y_n & = \Bigl\{
                \lambda;\ \lambda \in \N^n;\ 
                \biggl[\begin{aligned}
                \lambda_1 \geq \ldots \geq \lambda_d \geq 0\ \\
            \lambda_1 + \ldots + \lambda_d = n
                \end{aligned} \ \Bigr\}
            && \begin{aligned}
                    &\text{($n$-block Young diagrams} \\
                    &\quad\text{of height at most $d$)}
                    \end{aligned} \\
        Y^* & = \bigcup_{n \geq 2}  \{ \lambda;\ \lambda \in Y_n;\ 
                                \lambda_2 \neq 0 \} \\
        m_\lambda & = \prod_{1 \leq i < j \leq d}
            \frac{\lambda_i - \lambda_j + j - i}{j - i} 
            && \begin{aligned}
                &\text{(Dimension corresponding} \\
                & \quad \text{representation $\GL(d)$)}
            \end{aligned}\\
        c_k & 
         = \frac{1}{n} \sum_{l | n}
            d^\frac{n}{l}
       \mu\Bigl(\frac{l}{\gcd(l,k)}\Bigr)
                \frac{\phi(l)}{\phi\bigl(\frac{l}{\gcd(l,k)}\bigr)} 
                && \text{(Ramanujan sum)} \\
                \varphi & :\quad \text{Euler's totient} \\
                \mu & :\quad \text{M\"obius function} \\
            \end{align*}
        \end{minipage}}\end{center}
See appendix~\ref{apx:tables} for some decompositions computed using
these formulae.

\section{An Unordered Pair of Qubits}\label{s:upair-qubits}

The Hilbert space of a pair of
    qubits is~$\mathbb{C}^2 \otimes \mathbb{C}^2$.
Write~$\Hil=\mathbb{C}^2 \otimes \mathbb{C}^2$.
Let~$\sigma\colon \Hil \to \Hil$ denote the unitary map
that exchanges the two qubits:
\begin{equation*}
    \sigma \colon \ 
    \begin{aligned}
     &\ket{00} \mapsto \ket{00} &&
        &\ket{01} \mapsto \ket{10}\\
     &\ket{10} \mapsto \ket{01} &&
        &\ket{11} \mapsto \ket{11}
    \end{aligned}
\end{equation*}
An important category to study semantics of finite-dimensional
quantum computation is~$\starpu^\op$, which we will define in a moment.
It is important as every object corresponds to a type of
finite-dimensional quantum system and every arrow to a program of
the corresponding type.  Conversely, every physical finite dimensional
quantum type and  program corresponds to an object and arrow
(respectively) in this category.

The objects are finite-dimensional C$^*$-algebras\footnote{Which
are automatically unital.}. As the norm plays no
r\^ole in this paper, these are equivalently semisimple $*$-algebras
over $\C$, and are also equivalently finite dimensional W$^*$-algebras.
Accordingly, we refer to them as $*$-algebras in the rest of the
paper. We remark, however, that there are finite-dimensional
$*$-algebras in the axiomatic sense that are not C$^*$-algebras and
these are excluded.
The arrows in~$\starpu$ are completely positive unital linear maps,
and in $\starpu^\op$ they are in the opposite direction.
\footnote{The category~$\starpu^\op$ and its variations occur
under different names in the literature.
The category~$\cat{Star}_\mathrm{cP}^\op$
of finite dimensional~C$^*$-algebras with c.p.~maps in the opposite direction
is equivalent to
the category~$\mathbf{CP}^*[\mathbf{FHilb}]$ from \cite{bew103}.
If we restrict to subunital maps,
we call the category~$\cat{Star}_\mathrm{cPsU}^\op$,
which is equivalent to
the category~$\mathbf{Q}$
from~\cite{bew121}
and
the category~$\mathrm{CPM}_s$ from~\cite{bew143}.}

The~$*$-algebra of a pair of
    qubits is~$B(\Hil) \cong M_4$.
The map that exchanges the two qubits
    is given by:
\begin{equation*}
    B(\sigma) \colon  
        M_4 \leftarrow M_4, \quad
        a \mapsto \sigma^{-1}  a \sigma = \sigma a \sigma.
\end{equation*}
We claim the~$*$-algebra of an unordered pair of
        qubits must be the coequalizer of~$B(\sigma)$ and~$\id$.
This is the equalizer in~$\starpu$,
which is the following subalgebra of~$M_4$
\begin{equation*}
    E = \{ a; \ a \in M_4 ;\ \sigma a \sigma = a \} \subseteq M_4.
\end{equation*}
First note the analogy with the classical case:
to form an unordered pair of bits, one takes the quotient with
respect to the equivalence relation defined by permuting the bits,
which identifies $01$ and $10$. This is a coequalizer in the category
$\Sets$.
Why is the coequalizer used?  The definition gives the following
rule: for every program~$f \colon M_4 \leftarrow \mathscr{A}$
invariant under swapping ($\sigma \after f = f$)
there is a unique
lift~$f' \colon E \leftarrow \mathscr{A}$
such that~$e \after f' = f$, where~$e\colon E \subseteq M_4$
is the coequalizer map.

What~$*$-algebra is~$E$?
We write~$\Sil$ for the \emph{symmetric part} of~$\Hil$:
\begin{equation*}
    \Sil = \{ v; \ v \in \Hil;\ \sigma v = v\}.
\end{equation*}
One might expect $E = B(\Sil)$, but this is not the case. There is
another summand of~$E$.  First, we must take a small detour.
It is easy to verify that the projection onto~$\Sil$ is given by
\begin{equation*}
    P_\Sil \colon v \mapsto \frac{v + \sigma v}{2},
\end{equation*}
which is called the \emph{symmetrizer}.
The complementary projection~$P_\Ail = I-P_\Sil$
\begin{equation*}
    P_\Ail \colon v \mapsto \frac{v - \sigma v}{2}
\end{equation*}
projects onto the antisymmetric
subspace of~$\Hil$, which is given by
\begin{equation*}
    \Ail = \{ v ; \ v\in \Hil;\ \sigma v = -v\}. 
\end{equation*}
By considering the images of the standard basis vectors
under~$P_\Ail$ and~$P_\Sil$, it is easy to determine that
\begin{equation*}
    \{ \ket{00}, \ket{11} ,
    \frac{1}{\sqrt{2}} \ket{01} + \frac{1}{\sqrt{2}} \ket{10}  \}
    \quad\text{and}\quad \{ 
    \frac{1}{\sqrt{2}} \ket{01} - \frac{1}{\sqrt{2}} \ket{10}  \}
\end{equation*}
are orthonormal bases for~$\Sil$ and respectively $\Ail$.

There is a map~$i\colon B(\Sil)\oplus B(\Ail) \to B(\Sil \oplus \Ail) \cong B(\Hil)$,
given by
\begin{equation*}
    (s, a) \mapsto s \oplus a = \begin{pmatrix}
                                    s & 0 \\ 0 & a
                                \end{pmatrix}.
\end{equation*}
Its image, $\IM i$, is actually the equalizer~$E$.  We have to show both
inclusions.

First, suppose~$a \in \IM i$.
Then~$a = P_\Sil a P_\Sil + P_\Ail a P_\Ail$.
Note that~$\sigma P_\Sil = P_\Sil \sigma = P_\Sil$
and~$\sigma P_\Ail = P_\Ail \sigma = -P_\Ail$.  Thus:
\begin{equation*}
    \sigma a \sigma  = \sigma P_\Sil a P_\Sil \sigma +
                            \sigma P_\Ail a P_\Ail \sigma 
                     = P_\Sil a P_\Sil + P_\Ail a P_\Ail = a.
\end{equation*}
Hence~$a \in E$.

Conversely, suppose~$a \in E$.
First note that~$a = P_\Sil a P_\Sil + P_\Ail a P_\Ail + P_\Sil a
P_\Ail + P_\Ail a P_\Sil$.
Now since $\sigma a \sigma = a$, we have:
\begin{equation*}
     P_\Sil a P_\Sil + P_\Ail a P_\Ail + P_\Sil a P_\Ail + P_\Ail a P_\Sil =
         P_\Sil a P_\Sil + P_\Ail a P_\Ail - P_\Sil a P_\Ail - P_\Ail a P_\Sil.
\end{equation*}
Thus~$P_\Sil a P_\Ail = -P_\Ail a P_\Sil$.
Their images are orthogonal, hence~$P_\Sil a P_\Ail = P_\Ail a P_\Sil = 0$.
So~$a = P_\Sil a P_\Sil + P_\Ail a P_\Ail$,
and hence~$a \in \IM i$.

Thus~$E \cong B(\Sil) \oplus B(\Ail) \cong M_3 \oplus \mathbb{C}$.
At first one might be surprised that the antisymmetric vector~$\frac{1}{\sqrt{2}}\ket{01}
    - \frac{1}{\sqrt{2}}\ket{10}$ of~$\Hil$
    is a possible state of an unordered pair of qubits,
    since~$\sigma$ changes its sign.
The explanation is simple:
in~$*$-algebras, two states that differ only by global phase are identified.
Thus the antisymmetric vector is symmetric up to global phase~$-1$.

An astute reader might note that we have proven a bit more:
the~$*$-algebra associated
to an unordered pair of~$d$-level quantum systems
is given by~$B(\Ail) \oplus B(\Sil)$ as well,
where~$\Ail,\Sil \subseteq \mathbb{C}^d \otimes \mathbb{C}^d$
are defined similarly.

\section{Unordered Tuples}\label{s:utuples}
In the previous section, we have shown how to characterize
the~$*$-algebra for a pair of qubits.  In this section,
we will generalize to arbitrary tuples.
We define an unordered~$n$-tuple of~$d$-level quantum systems as follows.
Consider the Hilbert space~$(\mathbb{C}^d)^{\otimes n}$.
A permutation of~$n$ elements~$\pi \in S_n$ acts on it in an obvious way,
by permuting the basis vectors as follows:
\begin{equation}
\pi \colon \ \ket{i_1 i_2 \ldots i_n}
    \mapsto \ket{i_{\pi^{-1}(1)} \ldots i_{\pi^{-1}(n)}}. \label{eq:bigrepr}
\end{equation}
The equalizer of all~$\pi \in S_n$ in~$\starpu$
is the~$*$-algebra for unordered~$n$-tuples of~$d$-level quantum systems.
It is given by the following subalgebra of~$B((\mathbb{C}^d)^{\otimes n})$
\begin{equation*}
    E = \{ a; \ \pi^{-1} a \pi = a \text{ for all }\pi \in S_n \}
    \subseteq B((\mathbb{C}^d)^{\otimes n}).
\end{equation*}
The final result is:
\begin{equation*}
    E \cong \bigoplus_{\mathclap{\substack{
        \lambda_1 \geq \ldots \geq \lambda_d \geq 0 \\
        \lambda_1 + \ldots + \lambda_d = n \\
        \lambda_1, \ldots, \lambda_d \in \mathbb{N}}}}
        M_{m_\lambda}
    \quad\text{where}\quad
    m_\lambda = \prod_{1 \leq i < j \leq d}
            \frac{\lambda_i - \lambda_j + j - i}{j - i}.
\end{equation*}
To prove this, we will first review some of the basics of representation theory
of finite groups.
Then we will introduce Schur-Weyl duality to prove the result.

A representation of a group is a pair $(V,\rho)$, where $V$ is a
vector space and $\rho : G \rightarrow \GL(V)$ is a group homomorphism.
Often, one refers to the vector space~$V$
as the representation instead of the group homomorphism. When
considering the action of $g \in G$ on vectors $v \in V$ it is
common to leave out the~$\rho$
and write~$gv$ instead of~$\rho(g)v$.

We now give some examples of representations. The vector
space~$(\mathbb{C}^d)^{\otimes n}$ is a representation of~$S_n$,
by the action given in equation~\eqref{eq:bigrepr}.  Another one
is that for any group~$G$, we can consider~$\rho_\mathrm{trivial}\colon
G \to \GL(\C)$ given by~$\rho_\mathrm{trivial}(g) = I$.  This is
called the \emph{trivial representation}.

Given two representations~$\rho\colon G \to \GL(V)$
and~$\sigma \colon G \to \GL(W)$
a morphism~$f$ from~$\rho$ to~$\sigma$
is a linear map~$f\colon U \to V$
such that~$\sigma(g) f = f \rho(g)$ for every~$g \in G$.
That is: linear maps that commute with the group actions
of the representations.

We can relate morphisms of representations to the equalizer that
we want to calculate as follows.
\begin{align*}
    \repsn(\cdn, \cdn) 
            & = \{a ;\ a \in B(\cdn);\ 
                        \pi a = a \pi \text{ for all } \pi \in S_n \} \\
            & = \{a ;\ a \in B(\cdn);\ 
                        \pi^{-1} a \pi = a \text{ for all } \pi \in S_n \} \\
            & = E.
\end{align*}

Given two representations~$\rho\colon G \to \GL(V), \sigma: G \to \GL(W)$,
one can define the direct sum representation on~$V\oplus W$
by~$(\rho , \sigma)(g)(v, w) = (\rho(g)(v), \sigma(g)(w))$.
A representation is called indecomposable if it is not
the direct sum in this way of two other representations.

Given a representation on a vector space~$V$ and a subspace~$U$,
one calls~$U$ invariant (under~$G$) if
for every~$u \in U$ and~$g \in G$ we have~$g u \in U$.
A representation on~$V$ is called irreducible 
if the only invariant subspaces are~$\{ 0 \}$ and~$V$ itself. This
intentionally implies that the unique representation on the
zero-dimensional vector space is not irreducible, for the same
reason that $1$ is not prime and $\emptyset$ is not connected as a
topological space.

A slightly surprising, but welcome, theorem is that
a representation of a finite group is indecomposable if and
only if it is irreducible.
Furthermore, every representation is uniquely the direct sum of
irreducible representations (up to isomorphism).  See \cite[Proposition
1.5]{fh}.

Thus there are distinct irreducible representations~$U_\lambda$
and natural numbers~$m_\lambda$, called multiplicities, such
that~$\cdn \cong \bigoplus U_\lambda^{\oplus m_\lambda}$
and hence
\begin{equation*}
    E \cong \bigoplus_{\lambda,\mu}
        \repsn(U_\lambda^{m_\lambda}, U_\mu^{m_\mu}).
\end{equation*}

Now, given a morphism between representations,
it is easy to see that its kernel and image are invariant.
Thus, the only morphisms between irreducible representations
are invertible or zero maps.
This is the first part of Schur's lemma.
Consequently the maps between non-isomorphic irreducible representations
are $0$ and do not contribute to the direct sum, giving
\begin{equation*}
    E \cong \bigoplus_{\lambda}
        \repsn(U_\lambda^{m_\lambda}, U_\lambda^{m_\lambda}).
\end{equation*}
The second part of Schur's lemma is the following observation.
Suppose we have an endomorphism~$f$ of an irreducible representation $V$.
Since the base field $\C$ is algebraically closed, $f$ must have
an eigenvalue~$\lambda$, which is to say that $f - \lambda I$ has
non-trivial kernel.  The map $f - \lambda I$ is itself a morphism
of representations, and since $V$ is irreducible, $\ker(f - \lambda
I) = V$ and so $f-\lambda I = 0$.  That is to say: $f = \lambda I$.
Thus endomorphisms of irreducible representations
are scalar multiples of the identity.
We deduce
\begin{equation*}
    E \cong \bigoplus_{\lambda} M_{m_\lambda}.
\end{equation*}
Thus, if we know the irreducible representations of~$S_n$
and their multiplicities in~$\cdn$,
then we know~$E$.
Schur-Weyl duality solves this problem for us.
It gives a correspondence between
the irreducible representations of~$S_n$ in~$\cdn$
and of~$\GL(d)$ in~$\cdn$.
The space~$\cdn$ is a representation of~$\GL(d)$, via the following action
\begin{equation*}
    g v_1 \otimes \ldots \otimes v_n = (g v_1) \otimes \ldots \otimes (g v_d).
\end{equation*}
Schur-Weyl duality asserts
\begin{equation*}
    \cdn 
     \cong \bigoplus_{\mathclap{\substack{
        \lambda_1 \geq \ldots \geq \lambda_d \geq 0 \\
        \lambda_1 + \ldots + \lambda_d = n \\
        \lambda_1, \ldots, \lambda_d \in \mathbb{N}}}}
                U_\lambda \otimes V_\lambda 
     \equiv \bigoplus_{\mathclap{\substack{
        \lambda_1 \geq \ldots \geq \lambda_d \geq 0 \\
        \lambda_1 + \ldots + \lambda_d = n \\
        \lambda_1, \ldots, \lambda_d \in \mathbb{N}}}}
        U_\lambda^{\oplus \dim V_\lambda}
\end{equation*}
where~$U_\lambda$ are irreducible representations of~$S_n$
and~$V_\lambda$ are irreducible representations of~$\GL(d)$.
See~\cite[Exercise 6.30]{fh}.
Thus~$m_\lambda = \dim V_\lambda$.
Together with the duality statement, we are given explicit constructions
for~$U_\lambda$ and~$V_\lambda$.
See \cite[Theorem 4.3]{fh} and \cite[\S6.1]{fh}.
From this one can derive\cite[Theorem 6.3 (1)]{fh} that
\begin{equation*}
    \dim V_\lambda = 
    \prod_{1 \leq i < j \leq d}
    \frac{\lambda_i - \lambda_j +j - i}
        {j-i}.
\end{equation*}
In particular, in the case of unordered~$n$-tuples of qubits, we see
    $\dim V_\lambda = \lambda_i - \lambda_j + j - i$
and hence
\begin{equation*}
    E \cong \begin{cases}
        \bigoplus_{1 \leq i \leq \frac{n}{2}+1 } M_{2i-1} &
                n\text{ even} \\
        \bigoplus_{1 \leq i \leq \frac{n+1}{2} } M_{2i} &
                n\text{ odd}.
        \end{cases}
\end{equation*}

\section{A 3-cycle of Qubits}\label{s:3cycle}
Unordered tuples are defined by quotienting out
the action of the symmetric group.
Similarly, we can define other types
by quotienting out the action of a subgroup of the symmetric group.
The methods of the previous section can be adapted to
this situation as well. 
We will consider cycles, which are not as interesting a type
as unordered tuples, but they serve as an example easily related
to regular combinatorics.

A~$3$-cycle of qubits is given by the equalizer
\begin{equation*}
    E = \{ a; \ a \in M_8 ;\ \pi^{-1} a \pi = a;\ \pi \in C_3 \leq S_3 \}
                    \subseteq M_8.
\end{equation*}
The cyclic subgroup~$C_3$ of~$S_3$
contains~$\{ (), (1\ 2\ 3), (1\ 3\ 2) \}$.
We can use the same argument as before
to derive that~$E \cong \bigoplus_i M_{m_i}$,
where~$m_i$ is the multiplicity of the~$i$th irreducible representation
of~$C_3$ in~$\mathbb{C}^2 \otimes \mathbb{C}^2 \otimes \mathbb{C}^2$.
However, Schur-Weyl duality will not help this time.  We need to determine
the multiplicities~$m_i$ in another way.

To this end, we recall the theory of characters.
Given a representation~$\rho\colon G \to \GL(V)$.
For each~$g \in G$ we can consider the trace~$\Tr \rho(g)$.
This yields a map~$\chi_V=\Tr \after \rho \colon G \to \mathbb{C}$,
which is called the \emph{character} of~$\rho$.

By the cyclic property of the trace,
we have for any character~$\chi$
that~$\chi(h^{-1}gh) = \chi (g h h^{-1}) = \chi(g)$.
Thus on the same conjugacy class, a character will give the same value.
Such a function is called a \emph{class function}.

Using Schur's lemma one can work out that
\begin{equation*}
    \dim \repG(V,W) = \frac{1}{\#G} \sum_{g\in G}
                            \overline{\chi_V(g)} \chi_W(g)
                    = \begin{cases}
                        1 & V \cong W \\
                        0 & V \not\cong W.
                    \end{cases}
\end{equation*}
Also, using spectral decomposition, we can
derive~$\chi_{V\oplus W} = \chi_V + \chi_W$.
Thus, for two such class functions~$\alpha, \beta \colon G \to \mathbb{C}$,
one is lead to define
\begin{equation*}
    (\alpha, \beta) = \frac{1}{\#G} \sum_{g\in G} \overline{\alpha(g)} \beta(g).
\end{equation*}
This is an Hermitian inner product on the class functions.
In fact, with respect to this inner product
\begin{enumerate}
    \item the characters of irreducible representations are
            an orthonormal basis of the class functions;
    \item a representation~$V$ is irreducible if and
                only if~$(\chi_V, \chi_V)=1$;
    \item there are as many irreducible representations as
            conjugacy classes \emph{and}
    \item the multiplicity of~$V$ in~$W$ is~$(\chi_V, \chi_W)$.
\end{enumerate}
See \cite[\S2.2 and Proposition 2.30]{fh}.

Thus, to determine the multiplicities of the irreducible representations
of~$C_3$
in~$\mathbb{C}^2 \otimes \mathbb{C}^2 \otimes \mathbb{C}^2$,
it is sufficient to determine the character
of~$\mathbb{C}^2 \otimes \mathbb{C}^2 \otimes \mathbb{C}^2$
and the characters of the irreducible representations of~$C_3$.

We determine the irreducible representations of~$C_3$ as follows.
As~$C_3$ is Abelian, its conjugacy classes are trivial.
Write~$\pi$ for the generator of~$C_3$ such that~$C_3=\{1, \pi, \pi^2\}$.
Thus, we are looking for~$\#C_3=3$ irreducible representations.
The trivial representation maps every group element to the identity matrix.
It has character~$(1,1,1)$.
Then we have two~$1$-dimensional representations given
by~$\pi \mapsto (\omega)$
and~$\pi \mapsto (\omega^2)$,
where~$\omega = e^{\frac{2}{3}i \pi}$.
Using the inner product, we can compute that these are distinct
irreducible representations. We summarize these results in a character table:
\begin{center}\begin{tabular}{llll}
    $C_3 \leq S_3$& $1$ & $\pi$ & $\pi^2$ \\\hline
    trivial & $1$ & $1$& $1$ \\
    first & $1$ & $\omega$& $\omega^2$ \\
    second & $1$ & $\omega^2$& $\omega$ \\
\end{tabular}\end{center}

Now we compute the character~$\chi$
of~$\mathbb{C}^2 \otimes \mathbb{C}^2 \otimes \mathbb{C}^2$.
This is particularly easy because of
the way the action is defined: the value of the character on~$g$
is the number of basis vectors fixed by~$g$.
Thus:
\begin{center}\begin{tabular}{llll}
    $C_3 \leq S_3$& $1$ & $\pi$ & $\pi^2$ \\\hline
$\mathbb{C}^2 \otimes \mathbb{C}^2 \otimes \mathbb{C}^2$
      & $8$ & $2$& $2$
\end{tabular}\end{center}
We compute
\begin{align*}
    (\chi_\text{trivial}, \chi) & = 4 &
    (\chi_\text{first}, \chi) &= 2&
    (\chi_\text{second}, \chi) &= 2.
\end{align*}
Thus the~$*$-algebra for a~$3$-cycle of qubits is given
by~$E=M_4 \oplus M_2 \oplus M_2$.

\section{Cycles}\label{s:cycles}
Now we characterize arbitrary cycles:
a~$n$-cycle of~$d$-level quantum systems is given by the equalizer
\begin{equation*}
    E = \{a;\ a\in B(\cdn);\ g^{-1} a g = a;\ g \in C_n \leq S_n \} \subseteq
            B(\cdn).
\end{equation*}
First, we compute the irreducible representation of~$C_n$.
Let~$\pi \in C_n$ be such that~$C_n = \{ 1, \pi, \pi^2, \ldots, \pi^n\}$.
Note that by commutativity, the conjugacy classes are trivial.
For any~$0 \leq k \leq n$, define a 1-dimensional representation~$\rho_k$ by
\begin{equation*}
    \rho_k \colon C_n \to \GL(\C) \qquad
    \pi^i \mapsto (\omega ^{ki}),
\end{equation*}
where~$\omega = e^{2\pi i/n}$.
Note that~$\rho_0$ is the trivial representation.
Now, observe
\begin{equation*}
    (\rho_i, \rho_i) = \frac{1}{n} \sum_{0 \leq i < n} |\omega^{ik} |^2 = 1
\end{equation*}
and~$\Tr \rho_j(\pi) \neq \Tr \rho_i(\pi)$ whenever~$i \neq j$,
so these are~$k$ distinct irreducible representations.
The character table is given by
\begin{center}\begin{tabular}{llllcl}
    $C_n \leq S_m$& $1$ & $\pi$ & $\pi^2$ & \ldots & $\pi^{n-1}$ \\\hline
    $\rho_0$ & $1$ & $1$ & $1$ & \ldots & $1$ \\
    $\rho_1$ & $1$ & $\omega$ & $\omega^2$ & \ldots & $\omega^{n-1}$ \\
    $\rho_2$ & $1$ & $\omega^2$ & $\omega^4$ & \ldots & $\omega^{2(n-1)}$ \\
             &&\vdots &&&\vdots \\
    $\rho_{n-1}$ & $1$ & $\omega^{n-1}$ & $\omega^{2(n-1)}$
        & \ldots & $\omega^{(n-1)^2}$ \\
\end{tabular}\end{center}

Now we will compute that character~$\chi$
of the representation on~$\cdn$.
The value of~$\chi (\pi^i)$ is the number of basis vectors
that are fixed by~$\pi^i$.

All of the basis vectors are fixed by~$1=\pi^0$, so~$\chi(1)=d^n$.
The only basis vectors fixed by~$\pi$
are of the form~$\ket{vv\ldots v}$.
The general case is more subtle.
For instance, suppose~$n=4$ and~$d=2$. Then~$\ket{0101}$ is fixed by~$\pi^2$.

Given~$0 \leq i < n$.
If a basis vector~$\ket{v_1 \ldots v_n}$ is fixed by~$\pi^i$,
then we must have~$v_j = v_{\pi^i(j)} = v_{{\pi^{2i}}(j)} = \ldots$
for any~$0 \leq j < n$.
If~$i$ is coprime to~$n$,
then~$\{0,\pi^i(0), \pi^{2i}(0), \ldots \}$ (the orbit of
the subgroup generated by $\pi^i$)
ranges over all indices and thus the basis vector must be of the
form~$\ket{vv\ldots v}$.
If~$j$ is not coprime to~$n$,
then~$\{1,2,\ldots,n\}$ splits into several equally
    sized orbits.
The size of each of them is the order of~$\pi^i$,
which equals~$\frac{n}{\gcd(i,n)}$.
Thus the number of orbits is~$\gcd(i,n)$.
On each of the orbits, the basis vector has the same value, but is otherwise
unrestricted.  Thus there are~$d^{\gcd(i,n)}$ basis vectors fixed by~$\pi^i$.
Thus~$\chi(\pi^i) = d^{\gcd(i,n)}$.

Now, we will compute the multiplicity of the~$k$th irreducible representation
in~$\rho$, which is given by~$(\chi_k, \chi)$:
\begin{align*}
    (\chi_k, \chi)
        & = \frac{1}{n} \sum_{0 \leq j < n} \omega^{jk} d^{\gcd(j,n)} \\
        & = \frac{1}{n} \sum_{l | n}
                \sum_{\substack{1 \leq j \leq n \\
                \gcd(j,n)=l}} \omega^{jk} d^l\\
        & = \frac{1}{n} \sum_{l | n}
                d^l
                \sum_{\substack{1 \leq j \leq n \\
                \gcd(j,n)=l}} \omega^{jk}. \numberthis \label{eq:cycle-raw}
\end{align*}
As~$l$ divides~$j$, we may substitute~$jl$ for~$j$ and get:
\begin{equation*}
    (\chi_k, \chi)
         = \frac{1}{n} \sum_{l | n}
                d^l
                \sum_{\substack{1 \leq jl \leq n \\
                \gcd(jl,n)=l}} \omega^{jlk} 
         = \frac{1}{n} \sum_{l | n}
                d^l
                \sum_{\substack{1 \leq j \leq \frac{n}{l} \\
                \gcd(j,\frac{n}{l})=1}} \omega^{jlk}.
\end{equation*}
In \cite{r}, Ramanujan introduced (what are now called) \emph{Ramanujan sums}:
\begin{equation*}
    c_n(m) = \sum_{\substack{1 \leq h \leq n \\
        \gcd(h,n)=1}} e\Bigl(\frac{hm}{n}\Bigr),
\end{equation*}
where~$e(x) = e^{2 \pi i x}$.
Note that~$\omega^{jlk}=e\bigl(\frac{jlk}{n}\bigr)$. Consequently
\begin{equation*}
    (\chi_k, \chi)
         = \frac{1}{n} \sum_{l | n}
            d^l c_{\frac{n}{l}}(k) 
         = \frac{1}{n} \sum_{l | n}
            d^\frac{n}{l} c_{l}(k).
\end{equation*}
H\"older gave a simple expression for~$c_l(k)$,
see~\cite[Theorem 272]{hw}:
\begin{equation*}
    c_l(k) = \mu\Bigl(\frac{l}{\gcd(l,k)}\Bigr)
                \frac{\phi(l)}{\phi\bigl(\frac{l}{\gcd(l,k)}\bigr)},
\end{equation*}
where~$\mu$ is the M\"obius function and~$\phi$ is Euler's totient.
Therefore:
\begin{equation*}
    (\chi_k, \chi)
         = \frac{1}{n} \sum_{l | n}
            d^\frac{n}{l}
       \mu\Bigl(\frac{l}{\gcd(l,k)}\Bigr)
                \frac{\phi(l)}{\phi\bigl(\frac{l}{\gcd(l,k)}\bigr)}.
\end{equation*}
There are two cases of particular interest,
which can be proven directly from~\eqref{eq:cycle-raw}:
\begin{itemize}
\item
If~$n$ is a prime number, then:
\begin{align*}
    (\chi_k, \chi) =
            \begin{cases}
                \frac{d^n+(n-1)d}{n} & k=0 \\
                \frac{d^n-d}{n} & k>0.
            \end{cases}
\end{align*}

\item
The multiplicity corresponding to the trivial representation is
\begin{equation*}
    (\chi_0, \chi)
    = \frac{1}{n} \sum_{l | n}
            d^\frac{n}{l} \mu(1) \phi(l) 
        = \frac{1}{n} \sum_{l | n}
            d^l \phi\bigl(\frac{n}{l}\bigr).
\end{equation*}
This is MacMahon's formula for counting the number
of possible necklaces with~$n$ beads, where we may choose
from~$d$ different colors of beads. See~\cite[4.63]{gkp}.
\end{itemize}

\section{Unordered Words}\label{s:uwords}
Classically, a word is just a~$n$-tuple for some~$n$.
To work out what should be an unordered word,
we simply work out what is an unordered~$n$-tuple.
In the quantum analogue, such a reduction does not work.
Again, we need to tune our methods to work out a suitable equalizer.

The Hilbert space for quantum words over a~$d$-level quantum system is
the infinite dimensional Hilbert space
\begin{equation*}
    \Hil := \bigoplus_{n \in \N} (\C^d)^{\otimes n}.
\end{equation*}
Note that it only contains sequences that are square summable.
The corresponding von Neumann algebra
is the set of all bounded operators~$B(\Hil)$.

We will define an action~$\rho_\Hil$
of~$\prod_{n \in \N} S_n$ on~$\Hil$ as follows.
\begin{equation*}
    \rho_\Hil(\pi_1, \pi_2, \ldots)(\ket{i_1 \ldots i_m})
    = \ket{i_{\pi^{-1}_m(1)} \ldots i_{\pi^{-1}_m(m)}}
\end{equation*}
We wish to compute the equalizer of the actions, which is simply given by
\begin{align*}
    E & = \{a;\ a \in B(\Hil);\ \pi^{-1}a\pi =a \text{ for all } \pi \in
                \prod_{n \in \N} S_n \} \\
        & = \breppsn (\Hil, \Hil), \numberthis\label{eq:uwords1}
\end{align*}
where~$\breppsn(\Hil,\Hil)$ denotes the morphisms of representations
that are bounded (as linear maps between Hilbert spaces).

We cannot simply apply the same techniques as in Section~\ref{s:utuples}.
There are various difficulties.
First, $\Hil$ is infinite dimensional and
the group~$\prod_n S_n$ is not finite
so it does not follow from the theory we used previously that~$H$ splits into irreducible
representations of~$\prod_n S_n$.
Secondly, the infinite product~$\bigoplus$ is not a coproduct anymore.
We will work around these issues \emph{ad hoc}. It is possible to
give $\prod_n S_n$ a compact topology using Tychonoff's theorem and
use the representation theory of compact groups, but we do not
pursue that direction.

Let~$i_n\colon S_n \to \prod_{n \in \N} S_n$ denote the obvious inclusion
and~$p_n \colon B(H) \to B((\C^d)^{\otimes n})$ the obvious projection.
Then~$p_n \after \rho_H \after i_n$ is the action we considered
in~\eqref{eq:bigrepr}.
Recall that
\begin{equation*}
(\C^d)^{\otimes n} \cong \bigoplus_{\lambda \in Y_n}
U_\lambda^{\oplus m_\lambda}\quad \text{where}\quad
m_\lambda = \prod_{1 \leq i < j \leq d}
    \frac{\lambda_i - \lambda_j + j - i}{j - i}
\end{equation*}
and $U_\lambda$ are distinct irreducible representations for~$S_n$
indexed by
\begin{equation*}
    Y_n = \Bigl\{
                \lambda;\ \lambda \in \N^n;\ 
                \biggl[\begin{aligned}
                \lambda_1 \geq \ldots \geq \lambda_d \geq 0\ \\
            \lambda_1 + \ldots + \lambda_d = n
                \end{aligned} \ \Bigr\},
\end{equation*}
which are called~$n$-block Young diagrams of height at most~$d$.
The diagram~$\lambda$ is often depicted
as a row of~$\lambda_1$ blocks,
then a row of~$\lambda_2$ blocks beneath it and so on.
All blocks are left justified.
For instance,~$(4,2,0)$ is written
as {\tiny\yng(4,2)}.

Note that~$U_\lambda$ for any~$\lambda \in Y_n$
is an irreducible representation
for~$\prod_{n \in \N} S_n$ as well,
since~$S_m$ acts trivially on~$U_\lambda$ if~$m \neq n$.
However, not all~$U_\lambda$ are distinct.

For each~$n\in \N$,
there is the trivial representation of~$S_n$.
They correspond
to the Young diagrams of height~$1$
({\tiny\yng(1)}, {\tiny\yng(2)}, {\tiny\yng(3)}, \ldots).
They are all isomorphic as representations of~$\prod_{n \in \N} S_n$.
The representation isomorphism between any two, is the unique non-zero
map between the~$1$-dimensional subspaces.
We will show all other representations are distinct.

The \emph{kernel} of a representation $(V,\rho)$,
is the subgroup of elements that
map to the identity operator, equivalently the kernel of $\rho$ as a group homomorphism.
If two representations are isomorphic
then their kernels and dimensions are the same.

Given~$n, m \in \N$ and~$\lambda \in Y_n$ and~$\mu \in Y_m$
with~$\lambda \neq \mu$ such that, without loss of
generality, $U_\lambda$ is not a trivial representation. Suppose~$n
= m$ and~$U_\lambda$ is isomorphic to~$U_\mu$
as representation of~$\prod_{n \in \N} S_n$.
Then it is also isomorphic via the same isomorphism
as representation of~$S_n=S_m$, which is a contradiction.
Thus~$U_\lambda$ and~$U_\mu$ are distinct.

For the remaining case, suppose~$n \neq m$.
Because $U_\lambda$ is not a trivial representation,
there is an element~$\pi \in S_n$
that is not in its kernel.
If~$U_\mu$ is a trivial representation,
then~$U_\lambda$ and~$U_\mu$ must be distinct
as they have different kernels.
If~$U_\mu$ is not a trivial representation,
then there is an element~$\pi' \in S_m$
that is not in its kernel.
By definition of the action on~$\Hil$,
every element of~$S_n$ is in the kernel of~$U_\mu$.
Thus~$U_\lambda$ and~$U_\mu$ have different kernel.
Hence they are distinct.

We have a direct sum decomposition of~$\Hil$
into irreducible representations of~$\prod_{n \in \N} S_n$:
\begin{equation*}
    \Hil = \bigoplus_{n \in \N} (\C ^d) ^{\otimes n}
    \cong   U_{\mathrm{trivial}}^{\oplus \omega} \oplus
              \bigoplus_{\csubstack{n \in \N \\
                                    \lambda \in Y_n \\
                                    h(\lambda) \neq 1}}
                                    U_\lambda^{\oplus m_\lambda},
\end{equation*}
where~$U_\mathrm{trivial}
        = U_{\tiny\yng(1)}
        = U_{\tiny\yng(2)} = \ldots$ is the trivial representation.
Write
\begin{equation*}
    Y^* = \bigcup_{n \in \N} \{ \lambda;\ \lambda \in Y_n;\ h(\lambda)\neq 1 \}.
\end{equation*}
Now recall~\eqref{eq:uwords1}:
\begin{align*}
    E & = \breppsn(\Hil,\Hil) \\
      &\cong \breppsn(
     U_{\mathrm{trivial}}^{\oplus \omega} \oplus
              \bigoplus_{\csubstack{n \in \N \\
                                    \lambda \in Y_n \\
                                    h(\lambda) \neq 1}}
                                    U_\lambda^{\oplus m_\lambda},
     U_{\mathrm{trivial}}^{\oplus \omega} \oplus
              \bigoplus_{\csubstack{n \in \N \\
                                    \lambda \in Y_n \\
                                    h(\lambda) \neq 1}}
                                    U_\lambda^{\oplus m_\lambda}).
\end{align*}
Using Schur's lemma and the fact that~$\oplus$ is a biproduct, we derive
\begin{align*}
    E 
      &\cong 
    \breppsn(
     U_{\mathrm{trivial}}^{\oplus \omega},
     U_{\mathrm{trivial}}^{\oplus \omega}) \\
     &\qquad \oplus
    \breppsn(
                \bigoplus_{\lambda \in Y^*}
                                    U_\lambda^{\oplus m_\lambda},
                \bigoplus_{\lambda \in Y^*}
                                    U_\lambda^{\oplus m_\lambda}) \\
    & \cong B(\ell^2) \oplus 
    \breppsn(
                \bigoplus_{\lambda \in Y^*}
                                    U_\lambda^{\oplus m_\lambda},
                \bigoplus_{\lambda \in Y^*}
                                    U_\lambda^{\oplus m_\lambda}).
\end{align*}
We have to be a bit more careful for the right-hand summand,
since~$\bigoplus$ is not a countable biproduct.
\begin{alignat*}{2}
    &\breppsn(
                \bigoplus_{\lambda \in Y^*}
                                    U_\lambda^{\oplus m_\lambda},
                \bigoplus_{\lambda \in Y^*}
                U_\lambda^{\oplus m_\lambda})\\
    &\qquad = \biggl\{ (a_{\lambda \mu});\ 
                                \Biggl[
                    \begin{aligned}
                        &a_{\lambda\mu} \in
                    \breppsn (U_\lambda^{\oplus m_\lambda},
                    U_\mu^{\oplus m_\mu}); \\
                    &(a_{\lambda\mu}) \in B(\bigoplus_{\lambda \in Y^*}
                                        U_\lambda^{\oplus m_\lambda});\ 
                    \lambda, \mu \in Y^*
            \end{aligned}\ \biggr\} &\quad&\text{(dfn.)} \\
    &\qquad = \biggl\{ (a_{\lambda\lambda });\ 
                                \Biggl[
                    \begin{aligned}
                        &a_{\lambda\lambda} \in
                    \breppsn (U_\lambda^{\oplus m_\lambda},
                    U_\lambda^{\oplus m_\lambda}); \\
                    &(a_{\lambda\lambda}) \in B(\bigoplus_{\lambda \in Y^*}
                                        U_\lambda^{\oplus m_\lambda});\ 
                    \lambda \in Y^*
            \end{aligned}\ \biggr\} &&\text{(Schur's lemma)} \\
    &\qquad = \biggl\{ (a_{\lambda\lambda });\ 
                                \Biggl[
                    \begin{aligned}
                        &a_{\lambda\lambda} \in
                    \breppsn (U_\lambda^{\oplus m_\lambda},
                    U_\lambda^{\oplus m_\lambda}); \\
                    & \sup_\lambda \|a_{\lambda\lambda}\| < \infty; \ 
                    \lambda \in Y^*
            \end{aligned}\ \biggr\} &&\text{($*$, see below)}\\
    &\qquad \cong \biggl\{ (a_{\lambda});\ 
                                \Biggl[
                    \begin{aligned}
                        &a_{\lambda} \in
                    M_{m_\lambda} \\
                    & \sup_\lambda \|a_{\lambda}\| < \infty; \ 
                    \lambda \in Y^*
            \end{aligned}\ \biggr\} &&\text{(reindexing)}\\
            & \qquad
            \cong \prod_{\lambda \in Y^*} M_{m_\lambda}.
\end{alignat*}
Consequently
\begin{equation*}
    E \cong B(\ell^2) \oplus \prod_{\lambda \in Y^*} M_{m_\lambda}.
\end{equation*}
For step~$*$, note that the inclusion~$\subseteq$ is easy, and the
other inclusion is can be carefully checked using the definition
of the direct sum and noting the cross terms are zero. We also
emphasize that the infinite product should be interpreted for C$^*$
or W$^*$-algebras, with the norm bounded (the C$^*$-sum). This is, in general, a
strict subalgebra of the infinite product in $\C$-algebras or rings.

\subsection*{Acknowledgments}
We would like to thank Sam Staton for suggesting the problem.

The first author has been financially supported by the Netherlands
Organisation for Scientific Research (NWO) under TOP-GO grant no.
613.001.013 (The logic of composite quantum systems).

\nocite{*}
\bibliographystyle{eptcs}
\bibliography{main}

\appendix

\section{Computed decompositions}
\label{apx:tables}
For easy reference,
we have computed\footnote{%
The script used for the computation can
be found here: \url{https://westerbaan.name/~bas/math/bags.py}}.
 the decompositions into matrix algebras
of the C$^*$-algebras for unordered pairs, triples and quads for
various types.

\begin{table}[h]
    \begin{subtable}[b]{0.22\textwidth}
\begin{tabular}{rll}
d & \tiny\yng(2) & \tiny\yng(1,1) \\ \hline
    2 & $M_3$ & $\C$ \\
    3 & $M_6$ & $M_3$ \\
    4 & $M_{10}$ & $M_6$ \\
    5 & $M_{15}$ & $M_{10}$ \\
    6 & $M_{21}$ & $M_{15}$ \\
    7 & $M_{28}$ & $M_{21}$ \\
    8 & $M_{36}$ & $M_{28}$ \\
    9 & $M_{45}$ & $M_{36}$ \\
    10 & $M_{55}$ & $M_{45}$ \\
\end{tabular}
\caption{unordered pairs}
    \end{subtable}
    \begin{subtable}[b]{0.3\textwidth}
\begin{tabular}{rlll}
d & \tiny\yng(3) & \tiny\yng(2,1) & \tiny\yng(1,1,1)\\  \hline
    2  & $M_4$ & $M_2$ &\\
    3  & $M_{10}$ & $M_8$ & $\C$ \\
    4  & $M_{20}$ & $M_{20}$ & $M_4$ \\
    5  & $M_{35}$ & $M_{40}$ & $M_{10}$ \\
    6  & $M_{56}$ & $M_{70}$ & $M_{20}$ \\
    7  & $M_{84}$ & $M_{112}$ & $M_{35}$ \\
    8  & $M_{120}$ & $M_{168}$ & $M_{56}$ \\
    9  & $M_{165}$ & $M_{240}$ & $M_{84}$ \\
    10  & $M_{220}$ & $M_{330}$ & $M_{120}$ \\
\end{tabular}
\caption{unordered triples}
    \end{subtable}
    \begin{subtable}[b]{0.4\textwidth}
\begin{tabular}{rlllll}
d & \tiny\yng(4) & \tiny\yng(3,1) & \tiny\yng(2,2)
        & \tiny\yng(2,1,1) & \tiny\yng(1,1,1,1)\\  \hline
    2  & $M_5$ & $M_3$ & $\C$ &\\
    3  & $M_{15}$ & $M_{15}$ & $M_6$ & $M_3$ \\
    4  & $M_{35}$ & $M_{45}$ & $M_{20}$ & $M_{15}$ & $\C$\\
5 & $M_{70}$ & $M_{ 105}$ & $M_{ 50}$ & $M_{ 45}$ & $M_{ 5}$ \\
6 & $M_{126}$ & $M_{ 210}$ & $M_{ 105}$ & $M_{ 105}$ & $M_{ 15}$ \\
7 & $M_{210}$ & $M_{ 378}$ & $M_{ 196}$ & $M_{ 210}$ & $M_{ 35}$ \\
8 & $M_{330}$ & $M_{ 630}$ & $M_{ 336}$ & $M_{ 378}$ & $M_{ 70}$ \\
9 & $M_{495}$ & $M_{ 990}$ & $M_{ 540}$ & $M_{ 630}$ & $M_{ 126}$ \\
10 & $M_{715}$ & $M_{ 1485}$ & $M_{ 825}$ & $M_{ 990}$ & $M_{ 210}$ \\
\end{tabular}
\caption{unordered quads}
    \end{subtable}
    \caption{Decompositions into matrix algebras 
                    of unordered pairs, triples and quads of various
                        types.}
\end{table}

\end{document}